\documentclass[titlepage,12pt]{amsart}
\usepackage[stable]{footmisc}
\makeatletter
\def\myfnt{\ifx\protect\@typeset@protect\expandafter\footnote\else\expandafter\@gobble\fi}
\makeatother

\usepackage{amsfonts}
 \usepackage{amssymb}
 \usepackage{amsmath,amsxtra,amsthm}

\usepackage[usenames]{color}

\usepackage[dvips]{graphicx}

\usepackage[all,cmtip,matrix,arrow,curve]{xy}

\usepackage{mathrsfs}

\DeclareFontFamily{OT1}{pzc}{}
\DeclareFontShape{OT1}{pzc}{m}{it}{<-> s * [1.15] pzcmi7t}{}
\DeclareMathAlphabet{\mathpzc}{OT1}{pzc}{m}{it}

\DeclareSymbolFont{bbold}{U}{bbold}{m}{n}
\DeclareSymbolFontAlphabet{\mathbbold}{bbold}

\newtheorem{theorem}{Theorem}
\newtheorem{deff}{Definition}

\newtheorem{proposition}{Proposition}
\newtheorem{example}{Example}
\newtheorem{lemma}{Lemma}

\newtheorem{cor}{Corollary}

\newtheorem{rem}{Remark}
\renewcommand{\proof}{{\bf Proof.}~}

\newcommand{\bqa}{\begin{eqnarray}}
\newcommand\eqa {\end{eqnarray}}
\newcommand{\beq}{\begin{eqnarray}}
\newcommand{\beqn}{\begin{eqnarray}\nonumber}
\newcommand{\eeq}{\end{eqnarray}}
\newcommand{\be}{\begin{array}}
\newcommand{\ee}{\end{array}}

   \newcommand\vf\varphi

 \newcommand{\Hom}{\mathrm{Hom}}
 
 \newcommand{\Ker}{\mathrm{Ker}}

 \newcommand{\md}{\mathrm{d}}

 \newcommand{\cL}{\mathcal{L}}

 \newcommand{\R}{{\mathbb R}}

  \def\g{{\mathfrak g}}

\title[Universal Cartan-Lie algebroid]{Universal Cartan-Lie algebroid of an anchored bundle with connection and compatible geometries}%

\author{Alexei Kotov}
\address{Alexei Kotov: Faculty of Science, University of Hradec Kralove, Rokitanskeho 62, Hradec Kralove
50003, Czech Republic }
\email{oleksii.kotovATuhk.cz}

\author{Thomas Strobl}
\address{Thomas Strobl: Institut Camille Jordan, Universit\'e Claude Bernard Lyon 1, Universit\'e de Lyon,
43 boulevard du 11 novembre 1918, 69622 Villeurbanne Cedex, France}
\email{stroblATmath.univ-lyon1.fr}
\vspace{50mm}

\date{December, 2017}
\begin{document}
\vspace{50mm}
\begin{abstract} Consider an anchored bundle $(E,\rho)$, i.e.~a vector bundle $E\to M$ equipped with a bundle map $\rho \colon E \to TM$ covering the identity. M.~Kapranov showed in the context of Lie-Rinehard algebras that there exists an extension of this anchored bundle to an infinite rank universal free Lie algebroid $FR(E)\supset E$. We adapt his construction to the case of an anchored bundle  equipped with an arbitrary connection, $(E,\nabla)$, and show that it gives rise to a unique connection  $\tilde \nabla$ on  $FR(E)$ which is compatible with its Lie algebroid structure, thus turning  $(FR(E), \tilde \nabla)$ into a Cartan-Lie algebroid. Moreover, this construction is universal: any connection-preserving vector bundle morphism from $(E,\nabla)$ to a Cartan-Lie Algebroid $(A,\bar \nabla)$ factors through a unique Cartan-Lie algebroid morphism from $(FR(E), \tilde \nabla)$ to $(A,\bar \nabla)$.

\vskip 1mm\noindent Suppose that, in addition, $M$ is equipped with a geometrical structure defined by some tensor field $t$ which is compatible with $(E,\rho,\nabla)$ in the sense of being annihilated by a natural $E$-connection that one can associate to these data. For example, for a Riemannian base $(M,g)$ of an involutive anchored bundle $(E,\rho)$,  this condition implies that $M$ carries a Riemannian foliation.
%In general, the compatibility of a tensor $t$ with $(E,\rho,\nabla)$ implies its adequate invariance transversal to $\rho(E)$.
It is shown that every $E$-compatible tensor field $t$ becomes invariant with respect to the Lie algebroid representation associated canonically to the Cartan-Lie algebroid $(FR(E), \tilde \nabla)$.

%\vskip 1mm\noindent Suppose that, in addition, $M$ is equipped with a geometrical structure compatible with $(E,\rho,\nabla)$, i.e.~$M$ carries a tensor field $t$ annihilated by the $E$-connection ${}^\tau \nabla$, ${}^\tau \nabla t = 0$, where ${}^\tau \nabla_s v := [\rho(s),v] - \rho(\nabla_vs)$ for all $s \in \Gamma(E)$ and $v \in \Gamma(TM)$. Then $\tilde \nabla$ gives rise to a likewise FR(E)-connection ${}^\tau \tilde \nabla$, which is a representation of its Lie algebroid structure on $TM$ and yields ${}^\tau \tilde \nabla t =0$.

\vspace{20mm}
\noindent \keywordsname : Universal Lie algebra, Lie algebroids, Cartan connections, Riemannian foliations.
\vspace{20mm}

\tableofcontents
\end{abstract}
\setcounter{page}{0}

\maketitle
%IF ONE LEAVES THIS, THERE IS NO NUMBER ON THE PAGE OF THE INTRO...\thispagestyle{plain}

%\tableofcontents

\clearpage

\setcounter{page}{1}

\section{Introduction}

\vskip 3mm\noindent Every vector space gives naturally rise to a free infinite-dimensional Lie algebra. Applying the same strategy to an anchored vector bundle
$$
\xymatrix{ E \ar[dr]\ar[rr]^\rho && TM\ar[dl]\\
    & M &   }
%\xymatrix{
%E \ar[dd] \ar[rd]^{\rho} \\
%& TM \ar[ld] \\
%M          }
$$
needs some more care due to compatibility with the anchor map $\rho$, which, as a simple consequence of the Lie algebroid axioms, is required to become a morphism of the brackets. This implies in particular that in general the image of the anchor map will increase within this process.  It is shown by M.~Kapranov in \cite{Kapranov07} that any anchored module, a module over a commutative algebra
together with a morphism of modules with values in the module of derivations of the algebra,
gives rise in a canonical way to a free infinite-dimensional Lie-Rinehart algebra.\footnote{A Lie-Rinehart algebra is an algebraic counterpart of a Lie algebroid, cf.~\cite{Rinehart63}.}
A free Lie-Rinehart algebra admits a natural filtration the associated graded algebra to which is the free Lie algebra in the category of modules over the same algebra generated
by this module. We apply the construction of Kapranov to the category of smooth real manifolds---the original paper operates with Lie-Rinehart algebras over arbitrary ground fields---and call the resulting Lie algebroid $FR(E)\to M$.

\vskip 3mm\noindent The main purpose of this article is, however, to extend this relation between an anchored bundle $E$ and its free Lie algebroid $FR(E)$ to the lifting of particular additional structures from $E$ to $FR(E)$ such that appropriate compatibility conditions are satisfied. For the case of a vector bundle connection, e.g., there is no natural compatibility condition to be required if the vector bundle is merely an anchored bundle; however, if it is a Lie algebroid, this changes: let $A \to M$ be a Lie algebroid and $\nabla$ a connection on $A$. Any connection on $A$ gives rise to a splitting $\iota_{\nabla}$ of the natural projection map $J^1(A) \to A$, where $J^1(A)$ is the 1-jet bundle of sections of $A$. On the other hand, $J^1(A)$ carries a natural Lie algebroid structure itself, induced from the one on $A$. The compatibility consists of asking that $\iota_{\nabla} \colon A \to J^1(A)$ is a Lie algebroid morphism \cite{Blaom04}, in which case we call the connection a Cartan connection and the couple $(A,\nabla)$ a \emph{Cartan-Lie algebroid}. This compatibility condition can be re-expressed \cite{Kotov-Strobl_paper1} as the vanishing of the following tensor \cite{Mayer-Strobl09}
\begin{equation}
S:=2 \mathrm{Alt} \langle \rho, F_\nabla \rangle + \nabla \left( {}^AT \right) \:, \label{S}
\end{equation}
where  $F_\nabla \in \Gamma(A^* \otimes A \otimes \Lambda^2T^*M)$ is the curvature of $\nabla$, the anchor is considered as a section $\rho \in \Gamma(A^*\otimes TM)$, so that the contraction and skew-symmetrisation are defined in an obvious way, and ${}^AT$ is the $A$-torsion of the simple $A$-connection ${}^A\nabla$ on $A$ defined by ${}^A\nabla_s(s') := \nabla_{\rho(s)}s'$ for all $s,s' \in \Gamma(A)$.

\vskip 3mm\noindent Theorem \ref{FR_Cartan_extension}, proven in this paper, is a refinement of the above-mentioned result of Kapranov: given any anchored bundle $E$ equipped with an arbitrary  connection $\nabla$ there is a unique Cartan connection $\tilde \nabla$ on the corresponding free Lie algebroid $FR(E)$ which extends the one on $E \subset FR(E)$. It is interesting to see that it is precisely the compatibility condition $S=0$ which fixes the extension to all of $FR(E)$ uniquely. We call $(FR(E),\tilde \nabla)$ the \emph{free Cartan-Lie algebroid}
generated by the anchored bundle with connection $(E,\nabla)$. $(FR(E),\tilde \nabla)$ has a universality property, moreover,  which we will specify further below.
We mention as an aside that albeit we deal only with smooth manifolds in this paper, a purely algebraic version of this theorem in the spirit of \cite{Kapranov07} is quite obvious.

\vskip 3mm\noindent For an anchored bundle with connection $(E,\nabla)$ there is a natural compatibility with any tensor field $t$ defined over its base $M$: define the $E$-connection ${}^E \nabla$ when acting on vector fields  $v \in \Gamma(TM)$ by means of ${}^E \nabla_s v := [\rho(s),v] - \rho(\nabla_vs)$ for all $s \in \Gamma(E)$ and extend this canonically to all tensor fields over $M$. It is natural to ask that $t$ should be annihilated by this $E$-derivative:
\begin{equation}\label{tinv}
{}^E \nabla t = 0 \, .
\end{equation}
The meaning of this condition becomes clearer with an example: suppose the image of the anchor map is involutive,
$[\rho(\Gamma(E)),\rho(\Gamma(E))]\subset \rho(\Gamma(E))$, then $\rho(\Gamma(E))$ defines a singular foliation on $M$. If $M$ is equipped with a metric $g$ satisfying ${}^E \nabla g = 0$, then this singular foliation is Riemannian, and in particular transversally invariant with respect to the foliation.  Similar statements hold true for other geometrical structures defined by means of a tensor field $t$ satisfying Equation \eqref{tinv}. We note in parenthesis, if $E$ carries in addition a Lie algebroid structure and the connection $\nabla$ is compatible with it in the sense of $S=0$, then ${}^E \nabla$ as defined above provides an honest Lie algebroid representation on $TM$, $T^*M$, and its tensor powers, and the compatibility with $t$ then simply implies that this tensor is invariant under this canonical representation.\footnote{We refer the reader to \cite{Kotov-Strobl_paper1} for proofs and further details about the statements in this paragraph.}

\vskip 3mm\noindent The second result of the present paper, formulated in some generalisation in Proposition \ref{extension_of_compatibility} below, is that geometrical structures $t$ on $M$ which are compatible with the anchored bundle with connection $(E,\nabla)$ in the sense of Equation \eqref{tinv} remain compatible also with respect to $(FR(E),\tilde \nabla)$, i.e.~they become invariant with respect to the canonical representation of the universal free Cartan-Lie algebroid $(FR(E),\tilde \nabla)$ on $TM$. Moreover, this construction is universal: \emph{every} connection-preserving morphism from $(E,\nabla,t)$ to a $t$-compatible Cartan-Lie algebroid $(A,\nabla, t)$ such that the base map is the identity factors through this free Cartan Killing Lie algebroid.

\vspace{3mm} \noindent Although it would be desirable to find conditions under which the free (Cartan-)Lie algebroid admits a finite-dimensional reduction, for the moment we leave this problem open.  A necessary condition for such a reduction is that the---unmodified finite dimensional---base $M$ of $FR(E)$ carries a singular foliation: in the infinite rank setting, involutivity of the image of  $\rho_{FR(E)}$ is  not sufficient for its integrability. In addition, even if $FR(E)$ admits a finite rank reduction, there are in general  additional obstructions for the Cartan structure to reduce to the quotient Lie algebroid, since not every finite rank Lie algebroid even admits a compatible connection.

\vskip 5mm \noindent {\bf Acknowledgements.} The research of A.K. was supported by grant no. 18-00496S of the Czech Science Foundation. The research of T.S. was supported by the project MODFLAT of the European Research Council (ERC), the NCCR SwissMAP of the Swiss National Science Foundation, and, in particular, by Projeto P.V.E.~88881.030367/
2013-01 (CAPES/Brazil).

\section{Anchored bundles and free Cartan-Lie algebroids}\label{section:free_Cartan-Lie}

\vskip 3mm\noindent %$CLie(M)$, $Anch_{conn}(M)$, $Anch_{c}(M)$
Let us denote by $Anch_{c}(M)$ the category whose objects are anchored bundles with connections and morphisms are connection-preserving bundle morphisms
commuting with the anchor maps.
Let $CLie(M)$ be the category of Cartan-Lie algebroids over $M$.
Every Cartan-Lie algebroid is an anchored bundle and every connection-preserving Lie algebroid morphism is a morphism of the underlying anchored bundle structures, thus
there is a natural forgetful functor
\beq\label{Forgetful_functor} CLie(M)\to Anch_{c}(M)\,.
\eeq
\begin{theorem}\label{FR_Cartan_extension} The functor (\ref{Forgetful_functor})
admits a left-adjoint functor
\beq\label{FL_functor}
FR\colon Anch_{c}(M) \to CLie(M)
\eeq
whose value at an anchored bundle with connection $(E,\rho, \nabla)$ is a Lie algebroid $FR(E)$ together with
a Cartan connection and an embedding
 of anchored bundles $\imath\colon E\to FR(E)$, called the free Lie algebroid generated by $E$. Thus
we have a natural isomorphism
\beq\label{FR_isomorphisms}
\Hom_{CLie(M)} (FR (E),A)=\Hom_{Anch_{c}(M)} (E, A)\,
\eeq
for every Cartan-Lie algebroid $A$.
\end{theorem}

\noindent In other words, for every connection-preserving morphism of anchored bundles $\phi\colon E\to A$ there exists a unique
Cartan-Lie algebroid morphism $\tilde\phi\colon FR(E)\to A$ such that the following diagram is commutative:
\beqn
\xymatrix{   && FR (E)\ar@{.>}[d]^{\tilde\phi}\\
    E\ar@{^{(}->}[urr]^\imath \ar[rr]^{\phi} && A   }
\eqa

\noindent\proof The proof will consist of several consecutive steps. First, we construct the free almost Lie algebroid $FR^{alm}(E)$
generated by an anchored bundle $E$. This follows by the same method as in \cite{Kapranov07}. By an almost Lie algebroid we shall mean an anchored bundle equipped with a skew-symmetric bilinear operation satisfying the Leibniz rule with respect to the anchor map $\rho$, and $\rho$ is, moreover,  a morphism of brackets.

\vskip 3mm\noindent
We start with %$\mathcal{F\!L}^{alm}(E/\R)$, %$\cL\!ie(E)$
$FL^{alm}(E/\R)$, %$\cL\!ie(E)$
the free almost Lie $\R$-algebra generated by the real vector space of sections of $E$, along with $FL^{alm}(E)$, the bundle of free almost Lie algebras generated
by $E$ as a bundle over $M$. %, whose space of sections we denote by $\mathcal{F\!L}^{alm}(E)$.
Recall that an almost Lie algebra is like a Lie algebra except that the bracket operation does not
necessarily respect the Jacobi identity.
Both $FL^{alm}(E/\R)$ and $FL^{alm}(E)$ are naturally graded, such that
the degree $d$ factors of $FL^{alm}(E/\R)$ and $\Gamma\left(FL^{alm}(E)\right)$ are spanned by brackets involving exactly $d$ elements.
It is easily seen that all homogeneous factors of $\Gamma\left(FL^{alm}(E)\right)$ are finite-rank projective modules, hence we obtain the required grading
of $FL^{alm}(E)$ as a vector bundle whose degree $d$ factors are finite-dimensional vector bundles, so that the fiber at $x\in M$ is naturally isomorphic to the free almost Lie algebra generated by $E_x$.
The free almost Lie algebroid $FR^{alm}(E)$, which we now want to construct, is the union of an increasing sequence of anchored finite-rank bundles
\beq\label{FR_sequence}
 FR^{alm}_{\le 1}(E)\subset FR^{alm}_{\le 2}(E)\subset \ldots \,,
\eeq
which are defined inductively starting from
 $FR^{alm}_{\le 1}(E) =E$ as follows:
% , $q_1 = \mathrm{id}$.
suppose $FR^{alm}_{\le d}(E)$ is constructed as an anchored bundle with the anchor $\rho_d$
together with a surjective homomorphism of  real vector spaces
\beqn
 FL^{alm}_{\le d}(E/\R)\equiv \bigoplus\limits_{i=1}^d FL^{alm}_{i}(E/\R) \stackrel{q_d}{\longrightarrow} \Gamma (FR^{alm}_{\le d}(E))\,,
\eeq
where $q_1$ is the identity map. Then we define $FR^{alm}_{\le d+1}(E)$ as an anchored bundle whose space of sections is the quotient of $FL^{alm}_{d+1}(E/\R)$ by the following relations:
\beqn
 [s,r] &=& 0\,, \hspace{3mm} s\in \Gamma (E)\, , \; r\in\Ker (q_d) \,, \\ \nonumber
 [fs, s']- [s, fs'] &=& \rho (s)(f)s' - \rho_d (s')(f)s   \,, \hspace{3mm} s\in \Gamma (E)\, , s'\in\Gamma (FR^{alm}_{\le d}(E)) \, .
\eeq
The bracket on $FL^{alm}_{\le d}(E/\R)$ descends to a bracket on smooth sections of $FR^{alm}(E)$. The anchor map is uniquely determined by requiring the morphism property.
 The multiplication on smooth functions is given by the formula
\beqn
f[s, s']=[fs, s']+\rho_{d'}(s')(f)s = [s, fs']-\rho_d (s)(f)s' \,,
\eeq
where $s$ and $s'$ are arbitrary sections of $FR^{alm}_{\le d}(E)$ and $FR^{alm}_{\le d'}(E)$, respectively.
The filtration $\{FR^{alm}_{\le d}(E)\}$ makes $FR^{alm}(E)$ into a filtered almost Lie algebroid, and the associated graded almost Lie algebroid is isomorphic to $FL^{alm}(E)$ with trivial anchor map.

\vskip 3mm\noindent
It is worth mentioning that in the end the image of the anchor becomes involutive and although the original anchored bundle does not necessarily
carry a (singular) foliation, one now obtains an involutive (singular) tangent distribution which contains the image of $\rho$ of the original $E$. If it were just for this
involutive tangent distribution, it could have been obtained equally, and more easily, by completing $\rho(E)$ by means of iterated Lie brackets of $\rho(\Gamma(E))$.

\vskip 3mm\noindent Now we extend the connection on $E$ to the free almost Lie algebroid obtained above. For any $s,s'\in \Gamma (E)$
we claim
\beq\label{conn_on_brackets}
\nabla [s,s'] =\cL_{s}\left(\nabla s'\right) - \cL_{s'}\left(\nabla s\right) - \nabla_{\rho(\nabla s)} s' + \nabla_{\rho(\nabla s')}s\,,
\eeq
where $\cL_s$ acts on $\Gamma(T^*M \otimes E) \ni \omega \otimes s'$ by means of
\beq \label{cL}
\cL_s \left(\omega \otimes s'\right) :=  \cL_{\rho(s)}\omega \otimes s' + \omega \otimes [s,s'] \, ,
\eeq  where in the last equation $\cL$ denotes the standard Lie derivative.
The expression (\ref{conn_on_brackets}) is well-defined as from its definition it follows that
\beqn\nabla\left( [s,fs']-\rho(s)(f)s'-f[s,s'] \right)= 0\eeq for any smooth function $f$. This gives rise to the connection on
$FR^{alm}_{\le d}(E)$ for $d=2$. Now we proceed analogously by induction for all $d$. Let us notice that in each step we automatically obtain $S=0$, where $S$ is defined by the same formula
(\ref{S}) as for a Lie algebroid. This becomes obvious by rewriting \cite{Kotov-Strobl_paper1}  $S$ according to the following formula: $S(s,s')=\cL_{s}\left(\nabla s'\right) - \cL_{s'}\left(\nabla s\right) - \nabla_{\rho(\nabla s)} s' + \nabla_{\rho(\nabla s')}s - \nabla [s,s']$, cf.~also \cite{Blaom04}.

\vskip 3mm\noindent The last task is to define a Cartan structure on the free almost Lie algebroid and finally show that it descends to the associated free Lie algebroid.

\vskip 3mm\noindent Given an almost Lie algebroid $L$, there is a unique almost Lie algebroid structure on $J^1 (L)$
compatible with prolongations $j_1 \colon \Gamma(L) \to J^1 (L)$, namely which satisfies $\rho(j_1(s))=\rho(s)$ as well as
\begin{equation} \label{jets-bracket}
[j_1(s),j_1(s')] = j_1([s,s'])
\end{equation}  for every $s,s'\in\Gamma(L)$.\footnote{Alternatively, there is a (somewhat sophisticated) description of the canonical almost Lie algebroid structure on the bundle of $k-$jets of $L$ in terms of supergeometry, see Remark \ref{super_description}
below.} The notion of a Cartan connection along with the formula for the compatibility tensor \eqref{S} does not need the bracket to obey the Jacobi identity.
Thus starting with a connection $\nabla$ which satisfies $S=0$, we obtain a morphism of almost Lie algebroids $L\to J^1 (L)$ which is determined by the corresponding splitting $\iota_\nabla \colon L \to J^1 (L)$ defined by
\beq \label{iota}
\iota_\nabla = j_1(s) + \nabla(s) \, ;
\eeq
 In this equation, $\nabla(s)$ is viewed as a section of $J^1 (L)$ by means of Bott's exact sequence $0 \to T^*M \otimes L \to J^1 (L) \to L \to 0$ \cite{Bott}.

 \vskip 3mm\noindent The \emph{Jacobiator}, which is a measure of the failure of the bracket to satisfy the Jacobi identity, is defined as
 $\mathrm{Jac}\,(s_1,s_2,s_3)=[s_1, [s_2, s_3]]+\mathrm{cycl}(s_1,s_2,s_3)$ for every triple of sections %$s_1, s_2, s_3$
  of $L$. From the definition of an almost Lie algebroid, it follows that $\mathrm{Jac}$ is totally anti-symmetric and $C^\infty (M)$-linear in all its arguments.
  Since $\iota_\nabla$ is in particular a morphism of the brackets, one has
  \beq\label{sigma_vs_Jacobiator} \iota_\nabla\circ \mathrm{Jac}= \mathrm{Jac}\circ \iota_\nabla\otimes\iota_\nabla\otimes \iota_\nabla\,.
  \eeq

\begin{lemma}\label{lem:Jac_is_cov_const}
$\mathrm{Jac}$ is a covariantly constant map.
\end{lemma}
\noindent\proof  From (\ref{jets-bracket}) we obtain
  \beq \label{j1Jac}
  \mathrm{Jac}\,\left(j_1(s_1),j_1(s_2),j_1(s_3)\right)=j_1 \left( \mathrm{Jac}\,(s_1,s_2,s_3)\right)
  \eeq
for all $s_1, s_2, s_3\in\Gamma (L)$. After establishing that
\beq [j_1(s),\omega'\otimes s'] = \cL_s (\omega'\otimes s') \label{cross-bracket}
\eeq
and using that $\rho$ is a morphism of the brackets, one finds for all $\omega\in\Omega^1 (M)$
\beq\label{only_nonzero}
  \mathrm{Jac}\,\left(j_1(s_1),j_1(s_2),\omega\otimes s_3)\right)= \omega\otimes\mathrm{Jac}\,\left(s_1,s_2,s_3\right)\,.
\eeq
From  (\ref{cross-bracket}) and the fact that $[ \omega\otimes s,\omega'\otimes s']= \iota_{\rho(s')}\omega \left(\omega'\otimes s\right) - \iota_{\rho(s)}\omega' \left(\omega\otimes s\right)$,
we conclude that
for all sections $s_1,s_2,s_3$ and 1-forms $\omega_1, \omega_2, \omega_3$
\begin{eqnarray}
\mathrm{Jac}\,\left(j_1(s_1),\omega_2\otimes s_2,\omega_3\otimes s_3)\right)&=&0 \: ,
 \nonumber \\ \mathrm{Jac}\,\left(\omega_1\otimes s_1,\omega_2\otimes s_2,\omega_3 \otimes s_3\right)&=&0
\nonumber \: .
\end{eqnarray}
So we see that the only non-vanishing term of (\ref{sigma_vs_Jacobiator}) is of the form (\ref{only_nonzero}).
To be more explicit, let $b_i$ be a local basis of sections of $L$ and denote the corresponding connection coefficients
by $\omega_i^j$, $\nabla b_i = \omega_i^j \otimes b_j$. Then according to \eqref{sigma_vs_Jacobiator} and \eqref{iota},
\beqn \nabla \mathrm{Jac} (b_i,b_j,b_k) = \mathrm{Jac} (j_1(b_i),j_2(b_j),\nabla b_k)+ \mathrm{cycl}(ijk) = \omega_k^l \otimes  \mathrm{Jac} (b_i,b_j,b_l) + \mathrm{cycl}(ijk) \, .
\eeq
But this equation implies indeed $\nabla (\mathrm{Jac})=0$,
 or, equivalently, that the following diagram is commutative:
\beqn
\xymatrix{   \Gamma \left(\Lambda^3 (L)\right) \ar[d]_{\nabla}\ar[rr]^{\mathrm{Jac}}&& \Gamma (L)\ar[d]^{\nabla}\\
            \Gamma \left(T^*M\otimes\Lambda^3 (L)\right) \ar[rr]_{\mathrm{\mathrm{id} \otimes \mathrm{Jac}}}&& \Gamma \left(T^*M\otimes L\right) }
\eeq
Here the connection on $\Lambda^3 (L)$ is extended
by the Leibniz rule. This completes the proof of Lemma \ref{lem:Jac_is_cov_const}. $\square$

\begin{cor}\label{cor:quotient_CL}
Any Cartan connection on an almost Lie algebroid preserves the
Jacobi ideal of $L$, i.e. the ideal of sections generated by the image of the Jacobiator, thus it gives rise to a Cartan connection on the quotient Lie algebroid whenever it exists.
\end{cor}

\noindent Let us come back to the free almost Lie algebroid generated by $E$; it is obvious that the Jacobi ideal of $FR^{alm}(E)$, i.e.~the ideal generated by the Jacobiator, inherits the filtration
from the free almost Lie algebroid, such that the degree $d$ factors are finite-rank modules over the algebra of smooth functions.
Hence the quotient of $FR^{alm}(E)$ by the Jacobi ideal is a Lie algebroid,
which we denote by $FR(E)$. In \cite{Kapranov07} it is called the free Lie algebroid generated by $E$.
By Lemma \ref{lem:Jac_is_cov_const} and Corollary \ref{cor:quotient_CL}, we obtain the unique Cartan connection on $FR(E)$ which is compatible with the inclusion $E\hookrightarrow FR(E)$.
 We leave it to the reader to verify functorial properties of this construction.
$\blacksquare$

\begin{rem}\label{super_description}
Let $L$ be an almost Lie algebroid.
Consider $L[1]\to M$ as a graded superbundle with the degree $1$ odd fibers. In the way similar to the Lie algebroid case \cite{Vaitrob97}, an almost Lie algebroid structure
is in one-to-one correspondence with a degree $1$ vector field $Q$ on $L[1]$, defined by the Cartan's formula, such that $Q^2$ commutes
with all smooth functions on the base $M$. However, in contrast to Lie algebroids, the odd vector field $Q$ is not necessarily homological, as $Q^2 =0$
is equivalent to the Jacobi identity for the almost Lie algebroid structure on $L$. Now the canonical prolongation of $Q$ to the total space of $J^k (L[1])=J^k (L)[1]$ determines an almost Lie algebroid structure on the space of $k-$jets of $L$ compatible with the given one on $L$.
\end{rem}

\section{Compatible tensor fields on the base manifold}

\vskip 3mm\noindent For every anchored bundle $E$ equipped with a vector bundle connection $\nabla$, there is a natural compatibility of $(E,\rho,\nabla)$ with any tensor field $t$ defined over the base $M$. Define ${}^E{}\nabla$ by ${}^E{}\nabla_s v = \cL_{\rho(s)}v-\rho(\nabla_vs) $ for arbitrary sections $s\in\Gamma(E)$ when acting on vector fields $v$ and extend it by duality and the Leibniz rule to all tensors: so, e.g., for 1-forms $\omega$, one obtains $ {}^E{}\nabla_s \omega = \cL_{\rho(s)}\omega+\iota_{\rho(\nabla(s))} \omega$, where, by definition, $\iota_{\omega'\otimes v} \omega = ( \iota_v \omega)\omega' $. Then the compatibility of the tensor $t$ with the anchored bundle with connection is provided by the condition \eqref{tinv}. Somewhat more generally:\footnote{For a motivation and consequences of the compatibility condition see \cite{Kotov-Strobl_paper1}.}
\begin{deff} \label{Defff}
Let $(E,\rho)$ be an anchored bundle on $M$ and $\nabla^1, \ldots, \nabla^m$ be vector bundle connections on $E$  which give rise to ${}^E{}\nabla^1, \ldots, {}^E{}\nabla^m$
on tensors on $M$ as defined above for a single connection $\nabla$. A tensor $t$ on $M$ of  type $(r,s)$ with $r+s=m$  is called compatible with $(E, \nabla^1, \ldots, \nabla^m)$ if ${}^E{}\nabla^{comb}(t)=0$, where
 ${}^E{}\nabla^{comb} = {}^E{}\nabla^1 \otimes \mathrm{Id} \otimes \ldots \mathrm{Id}+ \ldots + \mathrm{Id} \otimes \ldots \mathrm{Id} \otimes  {}^E{}\nabla^m$.   \end{deff}

\begin{rem} \label{Rem}
In the particular case that $E$ is an almost Lie algebroid and $\nabla$  a Cartan connection on it---as defined in the previous subsection---any $E$-connection as above is an almost Lie algebroid representation. In other words, its $E$-curvature vanishes or, equivalently, $ {}^E{}\nabla_{[s,s']}= [{}^E{}\nabla_s, {}^E{}\nabla_{s'}]$. This then also applies to the combined $E$-covariant derivative ${}^E{}\nabla^{comb}$ above, provided, certainly, that  $\nabla^1, \ldots, \nabla^m$ are Cartan connections for the almost Lie algebroid structure $(E,\rho,[\cdot,\cdot])$. \end{rem}

\begin{proposition}\label{extension_of_compatibility}
Let $(E,\rho, \nabla^1, \ldots, \nabla^m)$ be an anchored bundle with $m$ connections and $t$ be a tensor field on the base $M$ of  type $(r,s)$ with $r+s=m$,
such that the compatibility condition ${}^E{}\nabla^{comb}(t)=0$ holds true. Let us extend $\nabla^1, \ldots, \nabla^m$ to the corresponding Cartan connections $\tilde\nabla^1, \ldots, \tilde\nabla^m$ on $FR(E)$.
Then $( FR(E), \tilde\nabla^1, \ldots, \tilde\nabla^m))$ is compatible with $t$, or, equivalently, $t$ is invariant with respect to the natural Lie algebroid representation of this $m$-fold free Cartan Lie algebroid
\end{proposition}
\noindent\proof Given that each $\tilde\nabla^i$ is a Cartan connection on $FR(E)$, i.e.~it respects the Lie algebroid structure in the sense of $S=0$, we can make use of the observation made in Remark \ref{Rem}:
it follows that ${}^{FR(E)}{}\tilde \nabla^i$ as well as ${}^{FR(E)}{}\tilde \nabla^{comb}$ are
 Lie algebroid representations of $FR(E)$.
Since the free Lie algebroid is generated by $E$, i.e.~the space of its sections is spanned over smooth functions by all multiple brackets of sections of $E$,
the identity ${}^E{}\nabla^{comb}_s(t)=0$ for all $s\in\Gamma (E)$ implies
${}^{FR(E)}{}\nabla^{comb}_\xi (t)=0$ for all sections $\xi$ of $FR(E)$, which proves the desired property.
$\blacksquare$

\vskip 3mm\noindent  We add two examples.

\begin{example}[cf.~\cite{Michor2008}, 3.29, example 2]
Consider the trivial rank two bundle $E$ with fiber generators $e_1,e_2$ over the base $\R^2$ parametrised by $x$ and $y$. Define the anchor map by means of $\rho(e_1)=\partial_x$ and $\rho(e_2)=\chi \partial_y$, where the function $\chi$ is identically zero for $x \leq 0$ and equal to $\exp(-1/x^2)$ for $x>0$. Since the repeated commutators of these vector fields produce higher and higher powers of $1/x$ in front of $\rho(e_2)$, this is not finitely generated. Moreover, we see that for non-positive $x$ the integral curves are straight horizontal lines, while for strictly positive $x$ there is only one two-dimensional leaf. Since these two foliations are glued together along $x=0$, however, one does not obtain a singular foliation in the sense of a partition of $M$, which is $\R^2$ here, into embedded submanifolds. Correspondingly, there can be also no finite rank quotient of $FR(E)$, since any finite rank Lie algebroid carries a singular foliation. The situation certainly does not improve, if additional structures are added, like a connection on $E$ and a compatible tensor on the base.
\end{example}

\begin{example}
Consider the trivial rank three bundle $E$ with fiber generators $e_1,e_2,e_3$ over the base $\R^3$ parametrised by $x,y,z$. Define $\rho$ by means of the following equations: $\rho(e_1)=\partial_x$, $\rho(e_2)=x\partial_y-y\partial_x$, and $\rho(e_3)= x \partial_z - z \partial_x$. As any trivial bundle, $E$ carries a canonical flat connection, $\nabla(e_i)=0$ for all $i=1,2,3$. If we equip the base with $g=\md x^2+ \md y^2+\md z^2$, $(E,\rho,\nabla,g)$ are compatible in the sense of Defintion \ref{Defff} with $\nabla^a := \nabla$ for $a=1,2$. As is evident from Equation \eqref{conn_on_brackets}, the extension of $\nabla$ to the connection $\tilde \nabla$ defined on all of the free Cartan-Lie algebroid $FR(E)$ is flat as well. And due to Proposition \ref{extension_of_compatibility}, the data $(FR(E),\tilde \nabla,g)$ are compatible, i.e.~they form an (infinite-rank) Killing Lie algebroid in the nomenclature of \cite{Kotov-Strobl_paper1}. This free Cartan-Lie algebroid does permit a finite quotient: factoring $FR(E)$ by the ideal generated by triple brackets, $FR(E)_{\geq 3}$, one obtains the action Killing Lie algebroid $(\g\times \R^3, \nabla^{can}, g)$, where $\g$ is the 6-dimensional isometry Lie algebra of $g$ and $\nabla^{can}$ is the canonical flat connection on the action Lie algebroid.
\end{example}

%\vskip 3mm\noindent  We conclude this section with a more intricate example: {\bf if you have such an example}

%\section{Examples and Open Problems}


\begin{thebibliography}{10}


\bibitem{Blaom04} A.D.~Blaom. \textit{Geometric structures as deformed infinitesimal symmetries.} Trans. Amer. Math.Soc. {\bf 358}, 3651--3671, 2006.
%arXiv:math/0404313

%\bibitem{Blaom05} A. D. Blaom. \textit{Lie algebroids and Cartan's method of equivalence.} Trans. Amer. Math.Soc. {\bf 364}, 3071--3135. 2012.
%arXiv:math/0509071



\bibitem{Bott} R.~Bott. \textit{Notes on the Spencer resolution.} Harvard University, Cambridge, Mass., 1963.



%\bibitem{Crainic-Fernandes03} M.~Crainic, R. L.~Fernandes. \textit{Integrability of Lie brackets.} Ann. of Math. {\bf 157(2)}, 575--620, 2003.

%\bibitem{daSilva-Weinstein99} A.C.~da Silva, A.~Weinstein. \textit{Geometric Models for Noncommutative Algebras.}
%Berkeley Mathematics Lecture Notes {\bf 10}, 184pp., 1999.


%\bibitem{delHoyo-Fernandes15} M.~del~Hoyo, R. L.~Fernandes. \textit{Riemannian metrics on Lie groupoids.}
%arXiv:1404.5989

%\bibitem{Fernandes02} R. L.~Fernandes. \textit{Lie algebroids, holonomy and characteristic classes.} Adv. Math., {\bf 170(1)}, 119--179, 2002.

%\bibitem{Stanciu} J.~Figueroa-O'Farrill, S.~Stanciu. \textit{Gauged Wess-Zumino terms and equivariant cohomology.}
% Phys. Lett. {\bf B341}, 153--159, 1994.

%\bibitem{Henneaux-Teitelboim} M.~Henneaux, C.~Teitelboim.
%\textit{Quantization of Gauge Systems.}
%Princeton University Press, 1992.


%\bibitem{Hitchin03} N.~Hitchin. \textit{Generalized Calabi-Yau manifolds.} Q. J. Math. {\bf 54 (3)}, 281--308, 2003.


%\bibitem{Hullequivariant}
%  C.~M.~Hull and B.~J.~Spence,
%  \textit{The Geometry of the gauged sigma model with Wess-Zumino term}.
%  Nucl.\ Phys.\ B {\bf 353} (1991) 379.



\bibitem{Kapranov07} M.~ Kapranov. \textit{Free Lie algebroids and the space of paths.} Selecta Mathematica {\bf 13}, 277--319, 2007.
%arXiv:math/0702584

%\bibitem{Kotov-Salnikov-Strobl14} A.~Kotov, V.~Salnikov and T.~Strobl. \textit{2d Gauge Theories and Generalized Geometry.}
%JHEP {\bf 21}, 2014.
%arXiv:1407.5439

%\bibitem{Kotov-Strobl07} A.~Kotov, T.~Strobl. \textit{Characteristic classes associated to
%Q-bundles.} International Journal of Geometric Methods in Modern Physics {\bf 12}, 2015.
%arXiv:0711.4106


\bibitem{Kotov-Strobl_paper1} A.~Kotov, T.~Strobl. \textit{Geometry on Lie algebroids: Compatible geometrical structures on the base.}

%\bibitem{Kotov-Strobl14} A.~Kotov, T.~Strobl. \textit{Gauging without Initial Symmetry.} J. Geom. Phys. {\bf 99}, 184-189, 2016.
%arXiv:1403.8119

%\bibitem{Kotov-Strobl_CYMH15} A.~Kotov, T.~Strobl.
%\textit{Curving Yang-Mills-Higgs Gauge Theories.}
%Phys.Rev.D {\bf 92} 085032, 2015.  %arXiv:1510.07654

%\bibitem{Mackenzie05} K.~Mackenzie. \textit{General Theory of Lie Groupoids and Lie Algebroids.} Cambridge University Press, 501 pp., 2005.

%\bibitem{LG-Lavau-Strobl} C.~Laurent-Gengoux, S.~Lavau, T.~Strobl. \textit{Lie infinity algebroid structures on singular foliations with resolution.}
%In preparation.

\bibitem{Mayer-Strobl09} C.~Mayer, T.~Strobl. \textit{Lie Algebroid Yang Mills with Matter Fields.} J. Geom. Phys. {\bf 59}, 1613--1623, 2009.
%arXiv:0908.3161

%\bibitem{Moerdijk-Mrcun03}  I.~Moerdijk, J.~Mrcun. \textit{Introduction to Foliations and Lie Groupoids.}
%Cambridge University Press, 184pp., 2003.

%\bibitem{Molino88} P.~Molino. \textit{Riemannian Foliations.} Progress in Mathematics {\bf 73}, Birkh\"{a}user, Boston, 344 pp., 1988.

%\bibitem{O'Neill66} B.~O'Neill. \textit{The fundamental equations of a submersion.} Michigan Math. J.
%{\bf 13(4)}, 459--469, 1966.

%\bibitem{Pradines67} J.~Pradines. \textit{Th\'{e}orie de Lie pour les groupoides
%diff\'{e}rentiables. Calcul diff\'{e}rentiel dans la cat\'{e}gorie
%des  groupoides
%infinit\'{e}simaux.} Comptes rendus Acad. Sci. Paris {\bf 264 A}, 245--248, 1967.

\bibitem{Michor2008} P.~Michor. \textit{Topics in Differential Geometry (Graduate Studies in Mathematics).}
Graduate Studies in Mathematics (Book 93), AMS, 494 pp., 2008.

\bibitem{Rinehart63} G.~Rinehart. \textit{Differential forms for general commutative algebras.} Trans. Amer. Math. Soc. {\bf 108}, 195--222, 1963.

%\bibitem{Salnikov-Strobl14} V.~Salnikov, T.~Strobl. \textit{Dirac Sigma Models from Gauging.}
%J. High Energy Phys. {\bf 110}, 2013.
% arXiv:1311.7116



%\bibitem{Strobl_AYM04} T.~Strobl. \textit{Algebroid Yang-Mills Theories.} Phys.Rev.Lett. {\bf 93} 211601, 2004.
%arXiv:hep-th/0406215

%\bibitem{Stroblinprep} T.~Strobl. \textit{From singular foliations to higher gauge theories.} In preparation.


\bibitem{Vaitrob97} A.~Vaintrob. \textit{Lie algebroids and homological vector fields.} Russian Math. Surveys {\bf 52}, 428--429, 1997.






\end{thebibliography}
\end{document}